\documentclass[11pt,twoside]{amsart}
\newtheorem{thm}{Theorem}[section]
\newtheorem{cor}[thm]{Corollary}
\newtheorem{claim}[thm]{Claim}

\newtheorem{rem}{\bf{Remark}}[section]

\numberwithin{equation}{section}
\def\pn{\par\noindent}

\newcommand{\kf}{\Lambda_k}
\newcommand{\lkf}{\Lambda^{l}_{k}}
\newcommand{\llc}{\Lambda^{l}_{3}}
\newcommand{\mad}{\mathrm{mad}}
\begin{document}

%--------------------------------------------------
%%Don not change any thing in this part
\centerline{ \scriptsize \it \textbf{\Small Article in Press}}
\leftline{ \scriptsize \it Bulletin of the Iranian Mathematical
Society  Vol. {\bf\rm XX} No. X {\rm(}201X{\rm)}, pp XX-XX.}

\vspace{1.3 cm}

%----------------------------------------------------------------------------
\title
{$k$-forested choosability of graphs with bounded maximum average degree}
\author{Xin Zhang, Guizhen Liu$^*$ and Jian-Liang Wu}

\thanks{{\scriptsize
This research is supported in part by National Natural Science Foundation of China (Grants No.\ 10871119, 10971121, 61070230) and Research Fund for the Doctoral Programs of Higher Education of China (Grant No.\ 200804220001).\\
\hskip -0.4 true cm MSC(2000): Primary: 05C15; Secondary: 05C10
%$\dag$Her research supported by a grant of the research council of the University of Tehran.\\
Received: 21 April 2010, Accepted: 24 September 2010\\
$*$Corresponding author
\newline\indent{\scriptsize $\copyright$ 2011 Iranian Mathematical
Society}}}

\maketitle

%-----------------------------------------------------------------------------
%This part will be filled in by BIMS
\begin{center}
Communicated by\;
\end{center}
%----------------------------------------------

\begin{abstract}  A proper vertex coloring of a simple graph is $k$-forested if the graph induced by the vertices of any two color classes is a forest with maximum degree less than $k$. A graph is $k$-forested $q$-choosable if for a given list of $q$ colors associated with each vertex $v$, there exists a $k$-forested coloring of $G$ such that each vertex receives a color from its own list. In this paper, we prove that the $k$-forested choosability of a graph with maximum degree $\Delta\geq k\geq 4$ is at most $\left\lceil\frac{\Delta}{k-1}\right\rceil+1$, $\left\lceil\frac{\Delta}{k-1}\right\rceil+2$ or $\left\lceil\frac{\Delta}{k-1}\right\rceil+3$ if its maximum average degree is less than $\frac{12}{5}$, $\frac{8}{3}$ or $3$, respectively.
\end{abstract}

\vskip 0.2 true cm

%-----------------------------------------------------------------------------

\pagestyle{myheadings}
\markboth{\rightline {\scriptsize  X. Zhang, G. Liu and J.-L. Wu}}
         {\leftline{\scriptsize $k$-forested choosability of graphs with bounded maximum average degree}}

\bigskip
\bigskip

%-----------------------------------------------------------------------------
%-----------------------------------------------------------------------------

\vskip 0.4 true cm

\section{\bf Introduction}

\vskip 0.4 true cm

In this paper, all graphs considered are finite, simple and undirected. We use $V(G)$, $E(G)$, $\delta(G)$ and
$\Delta(G)$ to denote the vertex set, the edge set, the minimum degree and the maximum degree of a graph $G$,
respectively. The $maximum$ $average$ $degree$ of $G$ is defined by $\mad(G)=\max\{2|E(H)|/|V(H)|,
H\subseteq G\}$. Any undefined notation follows that of Bondy and Murty \cite{Bondy}.

A proper vertex coloring of $G$ is called an $acyclic$ $coloring$ of $G$ if there are no bichromatic
cycles in $G$ under this coloring. The smallest number of colors such that $G$ has an acyclic coloring is called
the $acyclic$ $chromatic$ $number$ of $G$, denoted by $\chi_a(G)$. This concept was introduced by Gr\"{u}nbaum
\cite{Grunbaum}, and has been extensively studied in many papers. A coloring such that for every vertex $v\in
V(G)$ no color appears more than $k-1$ times in the neighborhood of $v$ is called a $k$-$frugal$ coloring. The
notation of $k$-frugality was introduced by Hind et al. in \cite{Hind}.

Yuster mixed these two notions (setting $k=3$) in \cite{Yuster} and first introduced the concept of $linear$ $coloring$, which is a proper coloring of $G$ such that the graph induced by the vertices of any two color classes is the union of vertex-disjoint paths. The $linear$ $chromatic$ $number$ $lc(G)$ of the graph $G$ is the smallest number $t$ such that $G$ has a linear $t$-coloring. Linear coloring was also investigated by Esperet, Montassier and Raspaud in \cite{Esperet}, and by Raspaud and Wang in \cite{Raspaud}. In \cite{Esperet}, the authors introduced a concept of $k$-$forested$ $coloring$ of a graph $G$, which is defined to be a proper vertex coloring of $G$ such that the union of any two color classes is a forest of maximum degree less than $k$. So a linear coloring is equivalent to a $3$-forested coloring. The $k$-$forested$ $chromatic$ $number$ of a graph $G$, denoted by $\kf(G)$, is the smallest number of colors appearing in a $k$-forested coloring of $G$. Note that $\kf(G)=\chi_a(G)$ for $k>\Delta(G)$. If $L$ is an assignment of a list $L(v)$ of colors to each vertex $v\in V(G)$, then $G$ is said to be $k$-$forested$ $L$-$colorable$ if it has a $k$-forested coloring where each vertex is colored with a color from its own list. We say $G$ is $k$-forested $q$-$choosable$ if $G$ is $k$-forested $L$-colorable whenever $|L(v)|=q$ for every vertex $v\in V(G)$. The $k$-$forested$ $choice$ $number$ $\lkf(G)$ is the smallest integer $q$ such that $G$ is $k$-forested $q$-choosable. When $k=3$, this is just equivalent to the $linear$ $choice$ $number$, which has been investigated by Esperet et al. for the graphs with bounded maximum average degree \cite{Esperet}. Their result is as follows.

\begin{thm}\label{llc} {\rm \cite{Esperet}} Let $G$ be a graph with maximum degree $\Delta$.

\noindent $(1)$ If $\Delta\geq 3$ and $\mad(G)<\frac{16}{7}$, then
$\llc(G)=\left\lceil\frac{\Delta}{2}\right\rceil+1$.

\noindent $(2)$ If $\mad(G)<\frac{5}{2}$, then $\llc(G)\leq \left\lceil\frac{\Delta}{2}\right\rceil+2$.

\noindent $(3)$ If $\mad(G)<\frac{8}{3}$, then $\llc(G)\leq \left\lceil\frac{\Delta}{2}\right\rceil+3$.
\end{thm}

\vskip 0.4 true cm

This paper is devoted to the following extensions of Theorem \ref{llc}.

\begin{thm} \label{mad2}Given a positive integer $M\geq k\geq 4$, let $G$ be a graph with maximum degree $\Delta\leq
M$.

\noindent $(1)$ If $\mad(G)<\frac{12}{5}$, then $\lkf(G)\leq\left\lceil\frac{M}{k-1}\right\rceil+1$.

\noindent $(2)$ If $\mad(G)<\frac{8}{3}$, then $\lkf(G)\leq \left\lceil\frac{M}{k-1}\right\rceil+2$.

\noindent $(3)$ If $\mad(G)<3$, then $\lkf(G)\leq \left\lceil\frac{M}{k-1}\right\rceil+3$.
\end{thm}

\vskip 0.4 true cm

By the definition of the $k$-forested choice number and $k$-forested chromatic number, one can easily say that
$\lkf(G)\geq \kf(G)\geq \left\lceil\frac{\Delta}{k-1}\right\rceil+1$ for every graph $G$ with maximum
degree $\Delta$. Now setting $M=\Delta$ in Theorem \ref{mad2}, we have the following theorem as a corollary.

\begin{thm} \label{mad1}Let $G$ be a graph with maximum degree $\Delta\geq k\geq 4$.

\noindent $(1)$ If $\mad(G)<\frac{12}{5}$, then $\lkf(G)=\left\lceil\frac{\Delta}{k-1}\right\rceil+1$.

\noindent $(2)$ If $\mad(G)<\frac{8}{3}$, then $\lkf(G)\leq \left\lceil\frac{\Delta}{k-1}\right\rceil+2$.

\noindent $(3)$ If $\mad(G)<3$, then $\lkf(G)\leq \left\lceil\frac{\Delta}{k-1}\right\rceil+3$.
\end{thm}

\vskip 0.4 true cm

Since every planar or projective-planar graph $G$ with girth $g(G)$ satisfies
$\mad(G)<\frac{2g(G)}{g(G)-2}$, we obtain the direct corollary from Theorem \ref{mad1}.

\begin{cor}\label{col}
Let $G$ be a planar or projective-planar graph with maximum degree $\Delta\geq k\geq 4$.

$(1)$ If $g(G)\geq 12$, then $\lkf(G)=\left\lceil\frac{\Delta}{k-1}\right\rceil+1$.

$(2)$ If $g(G)\geq 8$, then $\lkf(G)\leq \left\lceil\frac{\Delta}{k-1}\right\rceil+2$.

$(3)$ If $g(G)\geq 6$, then $\lkf(G)\leq \left\lceil\frac{\Delta}{k-1}\right\rceil+3$.
\end{cor}

\begin{rem}
In Theorems $\ref{mad2}$ and $\ref{mad1}$, we always respectively assume $M\geq k$
or $\Delta\geq k$. That is because once when we assume $M<k$ or $\Delta< k$, then $\lkf(G)=\chi^l_a(G)$ holds for any
graph $G$, where $\chi^l_a(G)$ denotes the acyclic choice number of $G$.
\end{rem}

\vskip 0.4 true cm

\section{\bf {\bf \em{\bf Proof of Theorem \ref{mad2}}}}

\vskip 0.4 true cm

In Claim \ref{cm} below, we will use ($p$) to denote the relevant part of Theorem \ref{mad2} ($p=1,2,3$). For brevity we will write $Q=\left\lceil\frac{M}{k-1}\right\rceil$ and $q=Q+p$, so that in part ($p$) we wish to prove that $\lkf(G)\leq q$. Note that, since $M\geq k$,
\begin{equation} \label{eq1}
    Q\geq 2~~and~~q=Q+p\geq p+2.
\end{equation}

Suppose that part ($p$) of Theorem \ref{mad2} is false. Let $G$ be a minimal counterexample to it; that is, every proper subgraph $H$ of $G$ is $k$-forested $q$-choosable but $G$ itself is not. (Here note that $\mad(H)\leq \mad(G)$ if $H$ is a subgraph of $G$.) Let $L$ be a list assignment of a list $L(v)$ of $q$ colors to each vertex $v\in V(G)$, such that $G$ has no $k$-forested $L$-coloring.

By the minimality of $G$, every proper subgraph $H$ of $G$ has a $k$-forested $L$-coloring. If $c$ is a $k$-forested $L$-coloring of a proper induced subgraph $H$ of $G$, and $v\in V(G)$, we use $c(N_G(v))$ to denote the set of colors used by $c$ on neighbors of $v$, and $C_{k-1}(v)$ to denote the set of colors that are each used by $c$ on exactly $k-1$ neighbors of $v$. Note that if $v$ has at least one neighbor that is uncolored, then
\begin{equation} \label{eq2}
|C_{k-1}(v)|\leq \left\lfloor\frac{d_G(v)-1}{k-1}\right\rfloor\leq \left\lfloor\frac{\Delta-1}{k-1}\right\rfloor\leq \left\lfloor\frac{M-1}{k-1}\right\rfloor= \left\lceil\frac{M}{k-1}\right\rceil-1=Q-1.
\end{equation}

\begin{claim}\label{cm}
$G$ does not contain any of the following configurations:

\noindent{\rm(C1)} a $1$-vertex;

\noindent{\rm(C2)} a $2$-vertex adjacent to a ($\leq p$)-vertex;

\noindent{\rm(C3)} if $p\leq k-2$, a $2$-vertex adjacent to a ($\leq p+1$)-vertex and a ($\leq 2p+1$)-vertex;

\noindent{\rm(C4)} if $p=3$, a $4$-vertex adjacent to three or more $2$-vertices;

\noindent{\rm(C5)} if $p=3$, a $5$-vertex adjacent to five $2$-vertices.
\end{claim}

\begin{rem}
In proving Claim $\ref{cm}$, we assume only that $k\geq 2$ in {\rm(C1)}, $k\geq p+1$ in {\rm(C2)}, $k\geq \max\{p+2,4\}$ in {\rm(C3)}, $k,p\geq 3$ in {\rm(C4)}, and $k\geq 4$, $p\geq 3$ in {\rm(C5)}. These conditions certainly hold if the conditions given in {\rm(C3)}--{\rm(C5)} hold and also $p\leq 3$ and $k\geq 4$, as stated in Theorem $\ref{mad2}$.
\end{rem}

\begin{rem}
In each part of the following proof, we first delete a set of vertices $\{x_1,\ldots,x_n\}$ from $G$ to obtain an induced subgraph $H$ that satisfies Theorem $\ref{mad2}$, and then extend the coloring $c$ of $H$ to each of $x_1,\ldots,x_n$ one by one. One should be careful
here to update the color set $C_{k-1}(\cdot)$ each time $c$ has been extended. For example, the color set $C_{k-1}(\cdot)$ in terms of the coloring $c$ of $H$ may be different from the one in terms of the coloring $c$ of $H+v_1$ after extending $c$ to $v_1$, but we still use the same notation for simplicity.
\end{rem}

\vskip 0.4 true cm

\pn{\bf Proof.} (C1) Suppose $G$ contains a 1-vertex $v$. Let $c$ be a $k$-forested $L$-coloring of $G-v$, which exists by the minimality of $G$. Denote the neighbor of $v$ by $u$, and define
$$F(v):=\{c(u)\}\cup C_{k-1}(u).$$
Then $|F(v)|\leq Q$ by (\ref{eq2}), and so $L(v)\setminus F(v)\neq \emptyset$ since $|L(v)|=q>Q$ by (\ref{eq1}). So we can color $v$ with a color in $L(v)\setminus F(v)$, and the coloring obtained is a $k$-forested $L$-coloring of $G$, which is a contradiction.

(C2) Suppose $G$ contains a 2-vertex $v$ which is adjacent to a ($\leq p$)-vertex $u$. Let the other neighbor of $v$ be $w$. In view of (C1) we may assume that $p\geq 2$. Let $c$ be a $k$-forested $L$-coloring of $G-v$. Note that $C_{k-1}(u)=\emptyset$, since $p-1<k-1$. Define
$$F(v):=\left\{
          \begin{array}{ll}
            \{c(u)\}\cup c(N_G(u))\cup C_{k-1}(w), & \hbox{if $c(u)=c(w)$;} \\
            \{c(u),c(w)\}\cup C_{k-1}(w), & \hbox{if $c(u)\neq c(w)$.}
          \end{array}
        \right.
$$
Then, by (\ref{eq2}), $|F(v)|\leq 1+(p-1)+(Q-1)<Q+p=q=|L(v)|$, and so we can color $v$ with a color in $L(v)\setminus F(v)$. This gives a $k$-forested $L$-coloring of $G$, which is a contradiction.

(C3) Suppose $G$ contains a 2-vertex $v$ which is adjacent to a ($\leq p+1$)-vertex $u$ and a ($\leq 2p+1$)-vertex $w$. Let $c$ be a $k$-forested $L$-coloring of $G-v$. Note that $C_{k-1}(u)=\emptyset$, since $p<k-1$ by hypothesis. Define
$$F(v):=\left\{
          \begin{array}{ll}
            \{c(u)\}\cup [c(N_G(u))\cap c(N_G(w))]\cup C_{k-1}(w), & \hbox{if $c(u)=c(w)$;} \\
            \{c(u),c(w)\}\cup C_{k-1}(w), & \hbox{if $c(u)\neq c(w)$.}
          \end{array}
        \right.
$$
Let $i=|c(N_G(u))\cap c(N_G(w))|\leq |c(N_G(u))|\leq p$. If $c(u)=c(w)$ then
$$|F(v)|\leq 1+i+\left\lfloor\frac{2p-i}{k-1}\right\rfloor\leq 1+p+\left\lfloor\frac{p}{k-1}\right\rfloor=1+p<q=|L(v)|$$
by (\ref{eq1}), since $p<k-1$. So suppose $c(u)\neq c(w)$. If $p=1$ then $|c(N_G(w))|\leq 2$ and so $C_{k-1}(w)=\emptyset$, since $k-1>2$; thus $|F(v)|\leq 2<3\leq |L(v)|$ by (\ref{eq1}). If $p\geq 2$, then $|F(v)|<|L(v)|$ by the same argument as in (C2). In every case we can color $v$ with a color from $L(v)\setminus F(v)$ to get a $k$-forested $L$-coloring of $G$, which is a contradiction.

(C4) Suppose $p=3$ and $G$ contains a 4-vertex $v$ which is adjacent to three 3-vertices $x,y,z$. Denote the other neighbors of $v,x,y,z$ by $w,x',y',z'$ respectively. Let $c$ be a $k$-forested $L$-coloring of $G-\{v,x,y,z\}$. Clearly $C_{k-1}(v)=\emptyset$. Give $z$ a color $c(z)\in L(z)\setminus F(z)$ where
$$F(z):=\{c(w),c(z')\}\cup C_{k-1}(z');$$
this is possible since $|L(z)|\geq Q+3$ by (\ref{eq1}), while $|C_{k-1}(z')|\leq Q-1$ by (\ref{eq2}). Next, noting that $v$ has colored neighbors $z,w$ where $c(z)\neq c(w)$, and $C_{k-1}(u)=\emptyset$ for all $u\in N_G(v)\setminus \{w\}$, give $v$ a color $c(v)\in L(v)\setminus F(v)$ where
$$F(v):=\{c(w),c(z),c(x')\}\cup C_{k-1}(w).$$
Then, noting that $|C_{k-1}(v)|=\left\lfloor\frac{1}{k-1}\right\rfloor=0$ since $c(z)\neq c(w)$, give $y$ a color from $L(y)\setminus F(y)$ where
$$F(y):=\left\{
          \begin{array}{ll}
            \{c(v),c(w),c(z)\}\cup C_{k-1}(y'), & \hbox{if $c(v)=c(y')$;} \\
            \{c(v),c(y')\}\cup C_{k-1}(y'), & \hbox{if $c(v)\neq c(y')$.}
          \end{array}
        \right.
$$
Finally, noting that $c(v)\neq c(x')$, give $x$ a color from $L(x)\setminus F(x)$ where
$$F(x):=\{c(v),c(x')\}\cup C_{k-1}(v)\cup C_{k-1}(x'),$$
which is possible since now $C_{k-1}(v)\leq \left\lfloor\frac{2}{k-1}\right\rfloor\leq 1$. This result is a $k$-forested coloring of $G$, a contradiction.

(C5) Suppose $p=3$ and $G$ contains a 5-vertex $v$ which is adjacent to five 2-vertices $x_1,\cdots,x_5$. Denote the other neighbor of $x_i$
 by $x'_i$ ($i=1,\cdots,5$). Let $c$ be a $k$-forested $L$-coloring of $G-\{v,x_1,x_2,x_3,x_4\}$. (In fact we do not need $d(x_5)=2$, only assuming $d(x_5)<k$ is enough so that when we prepare to color $v$, $C_{k-1}(x_5)=\emptyset$.) Give $x_1$ a color $c(x_1)\in L(x_1)\setminus F(x_1)$ where
$$F(x_1):=\{c(x'_1),c(x_5)\}\cup C_{k-1}(x'_1),$$
then give $v$ a color $c(v)\in L(v)\setminus F(v)$ where
$$F(v):=\{c(x_1),c(x_5),c(x'_2),c(x'_3)\},$$
which is possible since $|L(v)|\geq p+2=5$ by (\ref{eq1}). Now, noting that $c(x_1)\neq c(x_5)$ so that (even after $x_2$ is colored) $|C_{k-1}(v)|\leq \left\lfloor\frac{2}{k-1}\right\rfloor=0$, and $c(v)\not\in\{c(x'_2),c(x'_3)\}$, give $x_i$ a color from $L(x_i)\setminus F(x_i)$ where
$$F(x_i):=\{c(v),c(x'_i)\}\cup C_{k-1}(x'_i)~~~~~~(i=2,3).$$
Finally, give $x_4$ a color from $L(x_4)\setminus F(x_4)$ where
$$F(x_4):=\left\{
            \begin{array}{ll}
              \{c(v),c(x_1),c(x_5)\}\cup C_{k-1}(x'_4), & \hbox{if $c(v)=c(x'_4)$;} \\
              \{c(v),c(x'_4)\}\cup C_{k-1}(v)\cup C_{k-1}(x'_4), & \hbox{if $c(v)\neq c(x'_4)$,}
            \end{array}
          \right.
$$
which is possible since now $|C_{k-1}(v)|\leq \left\lfloor\frac{3}{k-1}\right\rfloor\leq 1$. This result is a $k$-forested $L$-coloring of $G$, a contradiction.

\vskip 0.4 true cm

In the next, we will complete the proof of each part of Theorem \ref{mad2} by a discharging procedure applying to the minimal counterexample $G$ to the theorem. We involve the same idea during each of the three proofs (assign each vertex $v\in V(G)$ an initial charge $w(v)=d(v)$) and the only differences are the definition of the discharging rules and the estimation on the final charge $w^*(v)$ of each vertex $v$ in $G$.

\vskip 0.4 true cm

\pn{\bf Proof of Theorem \ref{mad2}(1)}. We define discharging rules as follows.

\vspace{1mm} \textbf{R1.1}. Each $3$-vertex gives $\frac{1}{5}$ to each adjacent $2$-vertex;

\textbf{R1.2}. Each $\geq 4$-vertex gives $\frac{2}{5}$ to each adjacent $2$-vertex.

\vspace{1mm} Since the configuration (C1) in Claim \ref{cm} is forbidden in $G$, we assume that
$d(v)\geq 2$ for any vertex $v\in V(G)$. Suppose $d(v)=2$. If $v$ is adjacent to a $2$-vertex, then by the forbiddance of configuration (C3) in $G$, $v$ receives
$\frac{2}{5}$ from its another neighbor; if $v$ is not adjacent to any $2$-vertex, then $v$ also receives at least
$\frac{2}{5}$ from its neighbors. So $w^*(v)\geq w(v)+\frac{2}{5}=\frac{12}{5}$, since $v$ gives nothing. Assume
that $d(v)=3$. By R1.1, it gives out at most $\frac{3}{5}$. So $w^*(v)\geq
w(v)-\frac{3}{5}=3-\frac{3}{5}=\frac{12}{5}$. Assume that $d(v)=d\geq 4$. By R1.2, it gives out at most
$\frac{2d}{5}$. So $w^*(v)\geq w(v)-\frac{2d}{5}=d-\frac{2d}{5}=\frac{3d}{5}\geq\frac{12}{5}$. Thus
$w^*(v)\geq\frac{12}{5}$ for each vertex $v\in V(G)$, proving that
\begin{align*}
\mad(G)&\geq \frac{2|E(G)|}{|V(G)|}=\frac{\sum_{v\in V(G)}d(v)}{|V(G)|}
=\frac{\sum_{v\in V(G)}w(v)}{|V(G)|}\\&=\frac{\sum_{v\in
V(G)}w^*(v)}{|V(G)|}\geq\frac{12|V(G)|/5}{|V(G)|}=\frac{12}{5}.
\end{align*}
This contradiction proves Theorem \ref{mad2}(1).

\vskip 0.4 true cm

\pn{\bf Proof of Theorem \ref{mad2}(2)}. We define discharging rules as follows.

\vspace{1mm} \textbf{R2.1}. Each $3$-vertex gives $\frac{1}{9}$ to each adjacent $2$-vertex;

\textbf{R2.2}. Each $d$-vertex($4\leq d\leq 5$) gives $\frac{1}{3}$ to each adjacent $2$-vertex;

\textbf{R2.3}. Each $\geq 6$-vertex gives $\frac{5}{9}$ to each adjacent $2$-vertex.

\vspace{1mm} Similarly as above, we assume that
$d(v)\geq 2$ for any vertex $v\in V(G)$. Suppose $d(v)=2$. Then $v$ cannot be adjacent to any 2-vertex since (C2) can not appear in $G$ by Claim \ref{cm}. If $v$ is adjacent to a $3$-vertex, then by the forbiddance of configuration (C3) in $G$, another neighbor
of $v$ must be a $(\geq 6)$-vertex, so $v$ receives totally $\frac{1}{9}+\frac{5}{9}=\frac{2}{3}$ by R2.1 and R2.3. If $v$ is not
adjacent to any $3$-vertex, then by R2.2 and R2.3, $v$ receives at least $\frac{1}{3}+\frac{1}{3}=\frac{2}{3}$.
So $w^*(v)\geq w(v)+\frac{2}{3}=\frac{8}{3}$, since $v$ gives nothing. Suppose $d(v)=3$. Then $v$ gives out at most
$\frac{1}{3}$ by R2.1, so $w^*(v)\geq w(v)-\frac{1}{3}=\frac{8}{3}$. Similarly, we can prove that $w^*(v)\geq
\frac{8}{3}$ for any ($\geq 4$)-vertex. Thus $w^*(v)\geq\frac{8}{3}$ for each vertex $v\in V(G)$, proving that
$\mad(G)\geq \frac{8}{3}$. This contradiction completes the proof of Theorem \ref{mad2}(2).

\vskip 0.4 true cm

\pn{\bf Proof of Theorem \ref{mad2}(3)}. We define discharging rules as follows.

\vspace{1mm}\textbf{R3}. Each $\geq 4$-vertex gives $\frac{1}{2}$ to each adjacent $2$-vertex.

\vspace{1mm}  Similarly we first assume $d(v)\geq 2$ for any $v\in V(G)$. Suppose $d(v)=2$. Then the two neighbors of $v$ must be $(\geq 4)$-vertices since the configuration (C2) in Claim \ref{cm} is forbidden in $G$.
Thus, $v$ receives together $1$ from its neighbors but gives nothing by R3, which implies that $w^*(v)\geq
w(v)+1\geq 3$. Suppose $d(v)=3$. Note that $v$ receives and gives nothing by R3, so $w^*(v)=w(v)=3$. Suppose
$d(v)=4$. By the forbiddance of configuration (C4) in $G$, $v$ can be adjacent to at most two $2$-vertices, so it gives out at most $2\times
\frac{1}{2}=1$ by R3. This implies $w^*(v)\geq w(v)-1=3$. Suppose $d(v)=5$. Noting that the configuration (C5) can not occur in $G$, $v$ can be
adjacent to at most four $2$-vertices, so it gives out at most $4\times \frac{1}{2}=2$ by R3. This implies
$w^*(v)\geq w(v)-2= 3$. Suppose $d(v)=t\geq 6$. We have $w^*(v)\geq w(v)-\frac{1}{2}t=\frac{1}{2}t\geq 3$ by R3.
Thus $w^*(v)\geq 3$ for each vertex $v\in V(G)$, proving that $\mad(G)\geq3$. This contradiction completes the proof of Theorem \ref{mad2}.

%-----------------------------------------------------------------------------
%-----------------------------------------------------------------------------

\vskip 0.4 true cm

\begin{center}{\textbf{Acknowledgments}}
\end{center}

The authors wish to appreciate the anonymous referees sincerely for their very helpful comments.

%-----------------------------------------------------------------------------
%-----------------------------------------------------------------------------

\bigskip
\bigskip

%{\bf Received: Month xx, 200x}
{\footnotesize \pn{\bf X. Zhang}\; \\ {School of
Mathematics}, {Shandong University}, {Jinan 250100, China.}\\
{\tt Email:
sdu.zhang@yahoo.com.cn}\\

{\footnotesize \pn{\bf G. Liu}\; \\ {School of
Mathematics}, {Shandong University}, {Jinan 250100, China.}\\
{\tt Email:
gzliu@sdu.edu.cn}\\

{\footnotesize \pn{\bf J.-L. Wu}\; \\ {School of
Mathematics}, {Shandong University}, {Jinan 250100, China.}\\
{\tt Email:
jlwu@sdu.edu.cn}\\
\end{document}